\documentclass[11pt,a4paper]{article}
\date{}

\usepackage{latexsym}
\usepackage{amsfonts}
\usepackage{amssymb, amsfonts, amsmath}
\usepackage{amsbsy}

\newtheorem{thm}{Theorem}[section]
 
 \newtheorem{lem}[thm]{Lemma}

  \newtheorem{defn}[thm]{Definition}
 \newtheorem{rem}[thm]{Remark}
 
 \numberwithin{equation}{section}
\begin{document}

\author{Tulkin H. Rasulov}

\title{\sc Discrete Spectrum of a Model Operator Related to Three-Particle
Discrete Schr\"{o}dinger Operators} \maketitle

\begin{abstract}
A model operator $H_\mu,$ $\mu>0$ corresponding to a three-particle
discrete Schr\"{o}dinger operator on a lattice $ \mathbb{Z}^3$ is
considered. We study the case where the parameter function $w$ has a
special form with the non degenerate minimum at the $n,\,n>1$ points
of the six-dimensional torus $\mathbb{T}^6.$ If the associated
Friedrichs model has a zero energy resonance, then we prove that the
operator $H_\mu$ has infinitely many negative eigenvalues
accumulating at zero and we obtain an asymptotics for the number of
eigenvalues of $H_\mu$ lying below $z,$ $z<0$ as $z\to -0.$
\end{abstract}



\section{INTRODUCTION}

We are going to discuss the following remarkable phenomenon of the
spectral theory of the three-body Schr\"{o}dinger operators, known
as the Efimov effect: if a system of three particles interacting
through pair short-range potentials is such that none of the three
two-particle subsystems has bound states with negative energy, but
at least two of them have a zero energy resonance, then this
three-particle system has an infinite number of three-particle bound
states with negative energy accumulating at zero.

For the first time the Efimov effect has been discussed in
\cite{Efimov}. An independent proof on a physical level of rigor has
been also given in \cite{Amad-Nob} and then many works devoted to
this subject, see for example,
\cite{Dell-Fig-Teta,Ovch-Sig,Sob,Tam-1,Tam-2}. A rigorous
mathematical proof of the existence of Efimov's effect was
originally carried out in \cite{Yaf}.

Denote by $N(z)$ the number of eigenvalues of the Hamiltonian lying
below $z,\,z < 0.$ The growth of $N(z)$ has been studied in
\cite{AHW} for the symmetric case. Namely, the authors of \cite{AHW}
have first found (without proofs) the exponential asymptotics of
eigenvalues corresponding to spherically symmetric bound states.
This result is consistent with the lower bound $\lim\limits_{z \to
-0} \inf  N(z) |log|z||^{-1}>0$ established in \cite{Tam-1} without
any symmetry assumptions.

In \cite{Sob} the asymptotics of the form $N(z)\sim {\cal U}_0
|log|z||$ as $z\to -0$ for the number $N(z)$ of bound states of a
three-particle Schr\"{o}dinger operator below $z,$ $z<0$ was
obtained, where the coefficient ${\cal U}_0$ depends only on the
ratio of the masses of the particles.

Recently  in \cite{Wang} the existence of the Efimov effect for
$N$-body quantum systems with $N\geq 4$ has been proved and a lower
bound on the number of eigenvalues was given.

In \cite{Abd-Lak,Alb-Lak-Mum,L-2,L-1,L-M} the presence of Efimov's
effect for the three-particle discrete Schr\"{o}din\-ger operators
has been proved and in \cite{Abd-Lak,Alb-Lak-Mum} an asymptotics for
the number of eigenvalues similarly to \cite{Sob,Tam-2} was
obtained.

In the present paper, we study the model operator $H_\mu,$ $\mu>0$
corresponding to a three-particle discrete Schr\"{o}dinger operator
on a lattice $\mathbb{Z}^3.$ Here we are interested to discuss the
case where the parameter function $w$ has a special form with the
non degenerate minimum at the $n,\,n>1$ points of the
six-dimensional torus $\mathbb{T}^6.$ If the associated Friedrichs
model has a zero energy resonance, then we prove that the operator
$H_\mu$ has infinitely many negative eigenvalues accumulating at
zero (in the considering case zero is the bottom of the essential
spectrum of $H_\mu$). Moreover, we establish the asymptotic formula
$$
\lim\limits_{z\to -0}\frac{N_\mu(z)}{|log|z||}=\frac{{\bf n}
\gamma_0}{4\pi}
$$
for the number $N_\mu(z)$ of eigenvalues of $H_\mu$ lying below $z,$
$z<0.$ Here the number ${\bf n}\equiv{\bf n}(n),$ ${\bf n}>1$ is
defined in Remark \ref{about value of v} (see below) and the number
$\gamma_0$ is a unique positive solution of the equation
\begin{equation}\label{gamma0}
\gamma \sqrt{3}
\cos h\frac{\pi \gamma}{2}= 8 \sin h \frac{\pi \gamma}{6}.
\end{equation}

The asymptotics obtained in this paper can be considered as a
generalization of the asymptotics, which was obtained in
\cite{Abd-Lak,Alb-Lak-Mum,ALM-3,Sob,Tam-2}. In \cite{ALM-3} the non
symmetric version of the operator $H_\mu$ was considered and the
spectrum of this operator was analyzed for an arbitrary function $w$
with $n=1.$

The organization of the paper is as follows. In Section 2 the model
operator $H_\mu$ is introduced as a bounded self-adjoint operator
and the main result of the paper is formulated. In Section 3 some
spectral properties of the  associated Friedrichs model
$h_\mu(p),\,p\in (-\pi, \pi]^3$ are studied. In Section 4, we reduce
the eigenvalue problem by the principle of Birman-Schwinger. Section
5 is devoted to the prove of the main result of the paper.

\section{MODEL OPERATOR AND STATEMENT OF THE\\ MAIN RESULT}

Let us introduce some notations used in this work. Denote by
${\mathbb{T}}^3$ the three-dimensional torus, the cube
$(-\pi,\pi]^3$ with appropriately identified sides. The torus
${\mathbb{T}}^3$ will always be considered as an abelian group with
respect to the addition and multiplication by real numbers regarded
as operations on the three-dimentional space ${\mathbb{R}}^3$ modulo
$(2 \pi {\mathbb{Z}})^3.$ Let $({\mathbb{T}}^3)^2={\mathbb{T}}^3
\times {\mathbb{T}}^3$ be a Cartesian product, $
L_2({\mathbb{T}}^3)$ be the Hilbert space of square-integrable
(complex) functions defined on ${\mathbb{T}}^3$ and
$L_2^s(({\mathbb{T}}^3)^2)$ be the Hilbert space of
square-integrable symmetric (complex) functions defined on
$({\mathbb{T}}^3)^2.$

Let us consider a model operator $H_\mu$ acting on the Hilbert space
$ L_2^s(({\mathbb{T}}^3)^2)$ as
$$
H_\mu=H_{0} - \mu V_1 - \mu V_2,
$$
where
$$
(H_{0}f)(p,q)=w(p,q)f(p,q),
$$
$$
(V_1f)(p,q)= \varphi(p) \int\limits_{{\mathbb{T}}^3} \varphi(s)
f(s,q)ds,
$$
$$
(V_2f)(p,q)= \varphi(q) \int\limits_{{\mathbb{T}}^3} \varphi(s)
f(p,s)ds.
$$

Here $\mu$ is a positive real number, the function $\varphi(\cdot)$
is a real-valued analytic even function on ${\mathbb{T}}^3$ and the
function $w$ has form
\begin{equation*}\label{formulas of w}
w(p,q)=\varepsilon(p)+\varepsilon(p+q)+\varepsilon(q)
\end{equation*}
with
\begin{equation*}\label{varepsilon}
\varepsilon(p)=\sum_{j=1}^3 (1-cos \,m p^{(j)}),\,p=(p^{(1)},
p^{(2)}, p^{(3)})\in {\mathbb{T}}^3,
\end{equation*}
where $m$ is the positive integer number.

Under these assumptions the operator $H_\mu$ is bounded and
self-adjoint in $ L_2^s(({\mathbb{T}}^3)^2).$

To formulate the main result of the paper we introduce the
Friedrichs model $h_\mu(p),\,p\in {\mathbb{T}}^3,$ which acts in $
L_2({\mathbb{T}}^3)$ as
$$
h_\mu(p)=h_{0}(p)-\mu v,
$$
where
$$
(h_{0}(p)f_1)(q)=w(p,q)f(q),
$$
$$
(vf)(q)=\varphi(q) \int\limits_{{\mathbb{T}}^3} \varphi(s)f(s)ds.
$$

The perturbation $\mu v$ of the operator $h_{0}(p),\,p\in
{\mathbb{T}}^3$ is a self-adjoint operator of rank one. Therefore in
accordance with the invariance of the essential spectrum under
finite rank perturbations the essential spectrum
$\sigma_{ess}(h_\mu(p))$ of $h_\mu(p),\,p\in {\mathbb{T}}^3$ fills
the following interval on the real axis:
$$
\sigma_{ess}(h_\mu(p))=[m(p); M(p)],
$$
where the numbers $m(p)$ and $M(p)$ are defined by
$$
m(p)= \varepsilon(p)+2 \sum_{j=1}^3 (1- cos \frac{m
p^{(j)}}{2}),\,p=(p^{(1)},p^{(2)},p^{(3)})\in {\mathbb{T}}^3,
$$
$$
M(p)= \varepsilon(p)+2 \sum_{j=1}^3 (1+cos \frac{m
p^{(j)}}{2}),\,p=(p^{(1)},p^{(2)},p^{(3)})\in {\mathbb{T}}^3.
$$

The following Theorem \cite{ALM-3} describes the location of the
essential spectrum of $H_\mu.$
\begin{thm}\label{ess of H} For the essential spectrum
$\sigma_{ess}(H_\mu)$ of the operator $H_\mu$ the equality
$$
\sigma_{ess}(H_\mu)=\bigcup\limits_{p\in
{{\mathbb{T}}^3}}\sigma_{disc}(h_\mu(p))\cup [0; \frac{27}{2}]
$$
holds, where $\sigma_{disc}(h_\mu(p))$ is the discrete spectrum of
$h_\mu(p),\,p\in {\mathbb{T}}^3.$
\end{thm}

\begin{defn}
The set $\bigcup\limits_{p\in
{{\mathbb{T}}^3}}\sigma_{disc}(h_\mu(p))$ resp. $[0; \frac{27}{2}]$
is called two- resp. three-particle branch of the essential spectrum
$\sigma_{ess}(H_\mu)$ of the operator $H_\mu,$ which will be denoted
by $\sigma_{two}(H_\mu)$ resp. $\sigma_{three}(H_\mu).$
\end{defn}

Denote by $n\equiv n(m)$ the number of the all points of the form
$(p_i, q_j)\in ({\mathbb{T}}^3)^2$ with $p_i=(p_i^{(1)}, p_i^{(2)},
p_i^{(3)})$ and $q_j=(q_j^{(1)}, q_j^{(2)}, q_j^{(3)})$ such that
$p_i^{(k)}, q_j^{(k)} \in \{0, \pm \frac{2}{m} \pi; \pm \frac{4}{m}
\pi; \cdots; \pm \frac{m'}{m} \pi\},$ $k=1,2,3$ and $p_s\neq p_l,$
$q_s\neq q_l$ for $s \neq l,$ where
\begin{equation*}\label{the number m'}
m'=\left \lbrace
\begin{array}{ll}
m-2,\,\, \mbox{if the number}\,\, m \,\, \mbox{is even} \\
m-1,\,\, \mbox{if the number}\,\, m \,\, \mbox{is odd} .
\end{array} \right.
\end{equation*}

It is easy to check that the function $w(\cdot,\cdot)$ has the
non-degenerate minimum at that points $(p_i, q_j)\in
({\mathbb{T}}^3)^2$ and $n=(m'+1)^6.$

Now we additionally assume that $m\geq 3.$ Because, it is easy to
show that, if $m=1,2,$ then $n=1.$ In this paper we are interested
to study the case where $n>1.$

We denote that $\overline{1,n}=\{1,2,\cdots,n\}.$

\begin{rem}\label{about value of v}
In our analysis of the discrete spectrum of $H_\mu$ crucial role is
played by the zeroes of the function $\varphi(\cdot)$ at the points
$q_j \in {\mathbb{T}}^3,$ $j=\overline{1,\sqrt{n}}$ (see, for
example \cite{ALM-3}). Suppose that at only ${\bf n},$ $1<{\bf
n}\leq n$ points of the set $\{q_j\}_{j=1}^{\sqrt{n}} $ the value of
the function $\varphi(\cdot)$ is nonzero. We consider the set $
\{(p_{s_i}, q_{s_i})\in ({\mathbb{T}}^3)^2: i=\overline{1,n}\},$
where $s_i=\overline{1,n},$ such that $\varphi(q_{s_i}) \neq 0,$
$i=\overline{1, {\bf n}}$ and $\varphi(q_{s_i})= 0,$
$i=\overline{{\bf n}+1, n}.$ Throughout this paper we shall use this
notation without further comments.
\end{rem}

\begin{rem}
Note that the equality $h_\mu(p_{s_1})\equiv
h_\mu(p_{s_i}),\,i=\overline{2,n}$ holds.
\end{rem}

Let $C({\mathbb{T}}^3)$ (resp. $L_1({\mathbb{T}}^3)$) be the Banach
space of continuous (resp. integrable) functions on
${\mathbb{T}}^3.$

\begin{defn}\label{virtual level} The operator $h_\mu(p_{s_1})$ is said to have a zero
energy resonance if the number 1 is an eigenvalue of the integral
operator
$$
(G\psi_\alpha)(q)=\frac{\mu \varphi(q)}{2}
\int\limits_{{{\mathbb{T}}^3}}
\frac{\varphi(s)\psi(s)ds}{\varepsilon(s)},\quad \psi \in
C({\mathbb{T}}^3)
$$
and at least one (up to normalization constant) of the associated
eigenfunctions $\psi$ satisfies the condition $\psi(q_{s_i})\neq
0,\,i=\overline{1,{\bf n}}.$
\end{defn}

\begin{rem}\label{about value of v1}
We notice that if the operator $h_\mu(p_{s_1})$ has a zero energy
resonance, then the function
\begin{equation}\label{f_0 and f_1}
f(q)=\frac{\mu \varphi(q)}{2 \varepsilon(q)}\in L_1({\mathbb{T}}^3)
\setminus L_2({\mathbb{T}}^3),
\end{equation}
obeys the equation $h_\mu(p_{s_1})f=0$ (see Lemma \ref{virtual level
other}).
\end{rem}

Set
$$
\mu_0=2\left(\,\,\int\limits_{{\mathbb{T}}^3}
\frac{\varphi^2(s)ds}{\varepsilon(s)}\right)^{-1}.
$$

\begin{rem}
We remark that the operator $h_\mu(p_{s_1})$ has a zero energy
resonance if and only if $\mu=\mu_0$ (see Lemma \ref{zero energy
resonance}).
\end{rem}

Let us denote by $\tau_{ess}(H_\mu)$ the bottom of the essential
spectrum $\sigma_{ess}(H_\mu)$ of $H_\mu$ and by $N_\mu(z)$ the
number of eigenvalues of $H_\mu$ lying below $z,$
$z<\tau_{ess}(H_\mu).$

\begin{rem} We note that $\tau_{ess}(H_{\mu_0}) = 0$ (see Lemma \ref{about bottom of ess
spec}).
\end{rem}

The main result of this paper is the following

\begin{thm}\label{main result} The operator $H_{\mu_0}$ has an infinitely
many negative eigenvalues accumulating at zero and for the function
$N_{\mu_0}(\cdot)$ the relation
\begin{equation}\label{main asymp}
\lim\limits_{z\to -0}\frac{N_{\mu_0}(z)}{|log|z||}=\frac{{\bf
n}\gamma_0}{4\pi}
\end{equation}
holds, where the number ${\bf n}$ is defined in Remark \ref{about
value of v} and the number $\gamma_0$ is a positive solution of the
equation (\ref{gamma0}).
\end{thm}

\begin{rem} Clearly, the infinite cardinality of the negative discrete spectrum of
$H_{\mu_0}$ follows automatically from the positivity of the number
${\gamma}_0.$
\end{rem}

\begin{rem} We point out that the asymptotics (\ref{main asymp})
is new and similar asymptotics have not yet been obtained for the
three-particle Schr\"{o}dinger operators on $\mathbb{R}^3$ and
$\mathbb{Z}^3.$
\end{rem}

\section{SPECTRAL PROPERTIES OF THE OPERATOR $h_\mu(p)$}

In this section we study some spectral properties of the Friedrichs
model $h_\mu(p),\,p\in {\mathbb{T}}^3,$ which plays a crucial role
in our analysis of the discrete spectrum of the operator $H_\mu.$

Let $\mathbb{C}$ be the field of complex numbers. For any $p\in
{\mathbb{T}}^3$ we define an analytic function $ \Delta_\mu(p\,;
\cdot)$ (the Fredholm determinant associated with the operator
$h_\mu(p),\,p\in {\mathbb{T}}^3$) in ${\mathbb{C}}\setminus
\sigma_{ess} (h_\mu(p))$ by

$$
\Delta_\mu(p\,; z)=1-\mu \int\limits_{{\mathbb{T}}^3}\frac{
\varphi^2(q)dq}{w(p,q)-z}.
$$

The following statement (see \cite{ALM-3}) establishes a connection
between of eigenvalues of $h_\mu(p),\,p\in {\mathbb{T}}^3$ and
zeroes of the function $ \Delta_\mu(p\,; \cdot),\,p\in
{{\mathbb{T}}}^3.$

\begin{lem}\label{eigenvalue} For any $p\in {\mathbb{T}}^3$ the operator $h_\mu(p)$
has an eigenvalue $z\in{\mathbb{C}}\setminus \sigma_{ess}
(h_\mu(p))$ if and only if $ \Delta_\mu(p\,; z)=0.$
\end{lem}

Since the function $w(\cdot, \cdot)$ has the non-degenerate minimum
at the points $(p_{s_i}, q_{s_i})\in ({\mathbb{T}}^3)^2,$
$i=\overline{1,n}$ and the function $\varphi(\cdot)$ is an analytic
function on ${\mathbb{T}}^3,$ the integral
$$
\int\limits_{{\mathbb{T}}^3} \frac{\varphi^2(q)dq}{w(p, q)},\, p\in
{\mathbb{T}}^3
$$
is finite.

By Lebesgue's dominated convergence theorem and the equality
$\Delta_\mu(p_{s_i}\,; 0)=\Delta_\mu(p_{s_1}\,; 0),$
$i=\overline{2,n}$ it follows that
$$
\Delta_\mu(p_{s_1}\,; 0)=\lim\limits_{p\to p_{s_i}}\Delta_\mu(p\,;
0),\,i=\overline{1,n}.
$$

We remark that the following three statements, which are useful for
the proof of main result can be proven similarly to corresponding
statements of \cite{Abd-Lak,ALM-3} and hence here for completeness
we only reproduce these statements without proofs.

\begin{lem}\label{zero energy resonance}
The operator $h_\mu(p_{s_1})$ has a zero energy resonance if and
only if $\mu=\mu_0.$
\end{lem}

\begin{lem}\label{main decomp}
The following decomposition holds
$$
\Delta_{\mu_0}(p\,; z)=2 \pi^2 \mu_0 \sum\limits_{j=1}^{{\bf n}}
\varphi^2(q_{s_j})
\sqrt{\frac{3}{4}|p-p_{s_i}|^2-z}+O(|p-p_{s_i}|^2)+O(|z|)
$$
as $|p-p_{s_i}|\to 0,$ $i=\overline{1,{\bf n}}$ and $z\to -0.$
\end{lem}

Set
$$
U_\delta(p_0) = \{p \in {\mathbb{T}}^3 : |p-p_0| < \delta\},\quad
p_0 \in {\mathbb{T}}^3,\quad \delta>0.
$$

\begin{lem}\label{ineq for Delta}
There exist positive numbers $C_1,\,C_2,\,C_3$ and $\delta$ such
that
$$
C_1 |p-p_{s_i}|^2 \leq |\Delta_{\mu_0}(p\,; 0)|\leq C_2
|p-p_{s_i}|^2, \quad p\in U_\delta(p_{s_i}),\quad i=\overline{{\bf
n}+1, n};
$$
$$
|\Delta_{\mu_0}(p\,; 0)|\geq C_3,\quad p\in {\mathbb{T}}^3 \setminus
\bigcup\limits_{i=1}^{n} U_\delta(p_{s_i}).
$$
\end{lem}

From the representation
$$
w(p,q)=|p-p_{s_i}|^2+(p-p_{s_i}, q-q_{s_i})+|q-q_{s_i}|^2+
O(|p-p_{s_i}|^4)+O(|q-q_{s_i}|^4)
$$
as $|p-p_{s_i}|,\,|q-q_{s_i}|\to 0,\,i=\overline{1,n}$ it follows
the following

\begin{lem}\label{main inequality} There exist the numbers
$C_1, C_2, C_3>0$ and $\delta>0$ such that \\
1) $C_1(|p-p_{s_i}|^2+|q-q_{s_i}|^2)\leq w(p,q) \leq
C_2(|p-p_{s_i}|^2+|q-q_{s_i}|^2)$ for $(p,q)\in U_\delta(p_{s_i})
\times U_\delta(q_{s_i}),$
$i=\overline{1,n};$\\
2) $w(p,q) \geq C_3$ for all $p,q,$ which at least one of the
conditions $p \not\in \bigcup\limits_{i=1}^n U_\delta(p_{s_i})$ and
$q \not\in \bigcup\limits_{i=1}^n U_\delta(q_{s_i})$ is fulfilled.
\end{lem}

\begin{lem}\label{about bottom of ess spec}
The operator $h_{\mu_0}(p),\,p\in {\mathbb{T}}^3$ has no negative
eigenvalues.
\end{lem}

{\bf Proof.}
 First we show that for any $p\in
{\mathbb{T}}^3\setminus \{p_{s_1}, p_{s_2}, \cdots, p_{s_n}\}$ the
inequality $\Delta_{\mu_0}(p\,; 0)>\Delta_{\mu_0}(p_{s_1}\,; 0)$
holds. Denote
$$
\Lambda(p)=\int\limits_{\mathbb{T}^3}\frac{\varphi^2(q)dq}{w(p,q)}.
$$

Since the functions $\varphi(\cdot)$ and $w(\cdot, \cdot)$ are even,
the function $\Lambda(\cdot)$ is also even. Then
$$
\Lambda(p)-\Lambda(p_{s_1})=\frac{1}{4}\int\limits_{\mathbb{T}^3}
\frac{2w(p_{s_1}, q)-(w(p, q)+w(-p,
q))}{w(p,q)w(-p,q)w(p_{s_1},q)}[w(p,q)+w(-p,q)]^2 \varphi^2(q)dq-
$$
\begin{equation}\label{ineq for w_2}
-\frac{1}{4}\int\limits_{\mathbb{T}^3}
\frac{[w(p,q)+w(-p,q)]^2}{w(p,q)w(-p,q)w(p_{s_1},q)} \varphi^2(q)dq.
\end{equation}

By the equalities
$$
w(p_{s_1}, q)-\frac{w(p, q)+w(-p,
q)}{2}=\sum_{j=1}^3(cos\,mp^{(j)}-1)(1+cos\,mq^{(j)})
$$
and (\ref{ineq for w_2}) we have the inequality $
\Lambda(p)-\Lambda(p_{s_1})<0$ for any $p\in {\mathbb{T}}^3\setminus
\{p_{s_1}, p_{s_2}, \cdots, p_{s_n}\}.$

By the definition of $\mu_0$ we have $\Delta_{\mu_0}(p_{s_1}\,;
0)=0.$ Hence the inequality
$$
\Delta_{\mu_0}(p\,;
z)>\Delta_{\mu_0}(p_{s_1}\,; 0)=0
$$
holds for any $p\in {\mathbb{T}}^3$ and $z<0.$ By Lemma
\ref{eigenvalue} it means that, the operator $h_{\mu_0}(p),\,p\in
{\mathbb{T}}^3$ has no negative eigenvalues. $\Box$

\begin{lem}\label{virtual level other}
The function $f,$ which is defined by (\ref{f_0 and f_1}), obeys the
equation $h_{\mu_0}(p_{s_1})f=0.$
\end{lem}

{\bf Proof.} First we show that $f\in L_1({\mathbb{T}}^3)\setminus
L_2({{\mathbb{T}}}^3),$ that is,
$$
\int\limits_{{\mathbb{T}}^3}|f(q)|dq<\infty \quad \mbox{and} \quad
\int\limits_{{\mathbb{T}}^3}|f(q)|^2dq=\infty.
$$

From the definition of $\mu_0$ it follows that
$\Delta_{\mu_0}(p_{s_1}\,; 0)=0.$ By the construction of the set
$\{(p_{s_i}, q_{s_i})\in ({\mathbb{T}}^3)^2: i=\overline{1,n}\} $ we
have that $\varphi(q_{s_i}) \neq 0,$ $i=\overline{1, {\bf n}}$ and
$\varphi(q_{s_i})= 0,$ $i=\overline{{\bf n}+1, n}.$

Using these facts and the definition of the function
$\varepsilon(\cdot)$ we obtain that there exist the numbers
$C_1,\,C_2,\,C_3>0$ and $\delta>0$ such that
$$
C_1|q-q_{s_i}|^2 \leq \varepsilon(q) \leq C_2|q-q_{s_i}|^2,\quad
q\in U_\delta(q_{s_i}),\quad i=\overline{1,n},
$$
$$
\varepsilon(q)\geq C_3,\quad q\in {\mathbb{T}}^3 \setminus
\bigcup_{i=1}^n U_\delta(q_{s_i}),
$$
$$
|\varphi(q)| \geq C_3,\quad q\in  U_\delta(q_{s_i}), \quad
i=\overline{1,{\bf n}}
$$
and in the case where ${\bf n}<n$ we have that
$$
C_1|q-q_{s_i}|^2 \leq |\varphi(q)| \leq C_2|q-q_{s_i}|^2,\quad q\in
U_\delta(q_{s_i}),\quad i=\overline{{\bf n}+1,n}.
$$

Applying latter inequalities we obtain that
$$
\int\limits_{{\mathbb{T}}^3}|f(q)|dq \leq C_1 \sum_{j=1}^{\bf n}\,
\int\limits_{U_\delta(q_{n_j})} \frac{dq}{|q-q_{n_j}|^2}+C_2<
\infty,
$$
$$
\int\limits_{{\mathbb{T}}^3}|f(q)|^2dq \geq C_1 \sum_{j=1}^{\bf n}\,
\int\limits_{U_\delta(q_{n_j})} \frac{dq}{|q-q_{n_j}|^4}+C_2=\infty.
$$

It is easy to check that the function $f$ obeys the equation
$h_{\mu_0}(p_{s_1})f=0.$ $\Box$

\section{THE BIRMAN-SCHWINGER PRINCIPLE}

For a bounded self-adjoint operator $A,$ acting in Hilbert space
${\cal R},$ we define $d(\lambda, A)$ as
$$
d(\lambda, A)=\sup \{dim F: (Au,u)>\lambda,\, u\in F \subset {\cal
R},\,||u||=1\}.
$$

$d(\lambda, A)$ is equal to the infinity, if $\lambda$ is in the
essential spectrum and if $d(\lambda, A)$ is finite, it is equal to
the number of the eigenvalues of $A$ bigger than $\lambda.$

By the definition of $N_{\mu}(z)$ we have
$$
N_{\mu}(z)=d(-z, -H_{\mu}),\,-z>-\tau_{ess}(H_{\mu}).
$$

Since the function $\Delta_{\mu}(\cdot\,; \cdot)$ is positive on
${\mathbb{T}}^3 \times (-\infty, \tau_{ess}(H_{\mu})),$ the positive
square root of $\Delta_{\mu}(p\,; z)$ exists for any $p\in
{\mathbb{T}}^3$ and $z < \tau_{ess}(H_{\mu}).$

In our analysis of the spectrum of $H_{\mu}$ the crucial role is
played the compact integral operator $T_{\mu}(z),\,z <
\tau_{ess}(H_{\mu}),$ which acts in $L_2({\mathbb{T}}^3)$ with the
kernel
$$
\frac{\mu \, \varphi(p)\varphi(q)}{\sqrt{\Delta_\mu(p\,; z)}
\sqrt{\Delta_\mu(q\,; z)}(w(p,q)-z)} .
$$

The following lemma is a realization of the well known
Birman-Schwinger principle for the operator $H_\mu$ (see
\cite{Abd-Lak,Alb-Lak-Mum,ALM-3,L-M,Sob,Tam-2,Wang}).

\begin{lem}\label{Birm-Schw} For $z < \tau_{ess}(H_{\mu})$
the operator $T_{\mu}(z)$ is compact and continuous in $z$ and one
has
$$
N_{\mu}(z) = d(1, T_{\mu} (z)).
$$
\end{lem}

This lemma has been proven in \cite{ALM-3} for the non symmetric
case.

\section{THE PROOF OF THE MAIN RESULT}

In this section we shall derive the asymptotics (\ref{main asymp})
for the number $N_{\mu_0}(z)$ of eigenvalues of the operator
$H_{\mu_0}$ lying below $z,$ $z<0,$ that is, we shall prove Theorem
\ref{main result}.

We shall first establish the asymptotics for $d(1, T_{\mu_0}(z))$ as
$z\to -0.$ Then Theorem \ref{main result} will be deduced by a
perturbation argument based on the following lemma.

\begin{lem}\label{A0+A1} Let $A(z)= A_0(z)+A_1(z),$ where $A_0(z)$ (resp. $A_1(z)$)
is compact and continuous in $z < 0$ (resp. $z\leq 0$). Assume that
for some function $f(\cdot),\,f(z)\to 0,\, z \to 0$ one has
$$
\lim\limits_{z\to -0}f(z)\,d(\gamma, A_0(z)) = l(\gamma),
$$
and is continuous in $\gamma>0.$ Then the same limit exists for
$A(z)$ and
$$
\lim\limits_{z\to -0}f(z)\,d(\gamma, A(z)) = l(\gamma),
$$
\end{lem}

For the proof of Lemma \ref{A0+A1}, see Lemma 4.9 of \cite{Sob}.

Let $T(\delta; |z|)$ be the integral operator which acts in
$L_2({\mathbb{T}}^3)$ with the kernel
$$
\frac{1}{2\pi^2}\sum_{i=1}^{{\bf n}}
\frac{\chi_\delta(p-p_{s_i})\chi_\delta(q-q_{s_i})(\frac{3}{4}|p-p_{s_i}|^2+
|z|)^{-\frac{1}{4}}(\frac{3}{4}|q-q_{s_i}|^2+|z|)^{-\frac{1}{4}}}
{|p-p_{s_i}|^2+(p-p_{s_i}, q-q_{s_i})+|q-q_{s_i}|^2+|z|}.
$$

Here $\chi_\delta(\cdot)$ is the characteristic function of
$U_\delta(0).$

The following  lemma can be proven using Lemmas \ref{main decomp} --
\ref{main inequality}.

\begin{lem}\label{T(z) compakt} For any $z\leq 0$ and small $\delta>0$
the error $T_{\mu_0}(z)-T(\delta; |z|)$ is Hilbert-Schmidt operator
and is continuous in the uniform operator topology at the point $z =
0.$
\end{lem}

The space of the functions $f$ having support in
$\bigcup\limits_{i=1}^{{\bf n}} U_\delta(p_{s_i}),$ is an invariant
subspace for the operator $T(\delta; |z|).$ Let $T_0(\delta; |z|)$
be the restriction of the operator $T(\delta; |z|)$ to this
subspace, that is, the integral operator acting in
$L_2(\bigcup\limits_{i=1}^{{\bf n}} U_\delta(p_{s_i}))$ with the
kernel
$$
\frac{1}{2\pi^2}\sum_{i=1}^{\bf n}
\frac{(\frac{3}{4}|p-p_{s_i}|^2+|z|)^{-\frac{1}{4}}
(\frac{3}{4}|q-q_{s_i}|^2+|z|)^{-\frac{1}{4}}}{|p-p_{s_i}|^2+(p-p_{s_i},
q-q_{s_i})+|q-q_{s_i}|^2+|z|}.
$$

Denote by $diag\{A_1, A_2,\cdots, A_{\bf n}\}$ the ${\bf n} \times
{\bf n}$ diagonal matrix with operators $A_1, A_2,\cdots, A_{\bf n}$
as diagonal entries.

Since the space $L_2(\bigcup\limits_{i=1}^{\bf n}
U_\delta(p_{s_i}))$ is an isomorphous to
$\bigoplus\limits_{i=1}^{\bf n} L_2(U_\delta(p_{s_i})),$ the
operator $T_0(\delta; |z|)$ can be written as diagonal operator
$$
T_0(\delta; |z|)=diag\{T_0^{(1)}(\delta; |z|), T_0^{(2)}(\delta;
|z|),\cdots, T_0^{({\bf n})}(\delta; |z|)\},
$$
where $T_0^{(i)}(\delta; |z|),$ $i=\overline{1, {\bf n}}$ is the
integral operator acting in $\bigoplus\limits_{i=1}^{\bf n}
L_2(U_\delta(p_{s_i}))$ with the kernel
$$
\frac{1}{2\pi^2} \,
\frac{(\frac{3}{4}|p-p_{s_i}|^2+|z|)^{-\frac{1}{4}}
(\frac{3}{4}|q-q_{s_i}|^2+|z|)^{-\frac{1}{4}}}{|p-p_{s_i}|^2+(p-p_{s_i},
q-q_{s_i})+|q-q_{s_i}|^2+|z|}.
$$

One verifies that the operator $T_0(\delta; |z|)$ is unitary
equivalent to the operator $T_1(r),$ $r=|z|^{-\frac{1}{2}}$ acting
in $\bigoplus\limits_{i=1}^{\bf n} L_2(U_r(0))$ as
$$
T_1(r)=diag\{T_1^{(1)}(r), T_1^{(2)}(r),\cdots, T_1^{({\bf n})}(r)
\},
$$
where $T_1^{(i)}(r),$ $i=\overline{1, {\bf n}}$ is the integral
operator acting in the $L_2(U_r(0))$ with the kernel
$$
\frac{1}{2\pi^2} \frac{1}{(\frac{3}{4}|p|^2+1)^{\frac{1}{4}}
(\frac{3}{4}|q|^2+1)^{\frac{1}{4}}(|p|^2+(p, q)+|q|^2+1)}.
$$

We note that the equivalence of these operators is performed by the
unitary dilation
$$
B_r=diag\{B_r^{(1)}, B_r^{(2)},\cdots, B_r^{({\bf n})}\}:
\bigoplus\limits_{i=1}^{{\bf n}} L_2(U_\delta(p_{s_i})) \to
\bigoplus\limits_{i=1}^{{\bf n}} L_2(U_r(0)).
$$

Here the operator $B_r^{(i)}: L_2(U_\delta(p_{s_i})) \to
L_2(U_r(0)),$ $i=\overline{1, {\bf n}}$ acting by
$$
(B_r^{(i)}f)(p)=r^{-\frac{3}{2}}f(\frac{1}{r}(p-p_{s_i})).
$$

Since the space $\bigoplus\limits_{i=1}^{{\bf n}} L_2(U_r(0))$ is an
isomorphous to $L_2(U_r(0)),$  we rewrite the operator $T_1(r)$ as
integral operator acting in $L_2(U_r(0))$ with the kernel
$$
\frac{{\bf n}}{2\pi^2} \frac{1}{(\frac{3}{4}|p|^2+1)^{\frac{1}{4}}
(\frac{3}{4}|q|^2+1)^{\frac{1}{4}}(|p|^2+(p, q)+|q|^2+1)}.
$$

Further, we may replace $(\frac{3}{4}|p|^2+1)^{\frac{1}{4}},$
$(\frac{3}{4}|q|^2+1)^{\frac{1}{4}}$ and $|p|^2+(p, q)+|q|^2+1$ by
$(\frac{3}{4}|p|^2)^{\frac{1}{4}}(1-\chi_1(p)),$
$(\frac{3}{4}|q|^2)^{\frac{1}{4}}(1-\chi_1(q))$ and $|p|^2+(p,
q)+|q|^2,$ respectively, we have the operator $T_2(r).$ The error
$T_1(r)-T_2(r)$ will be a Hilbert-Schmidt operator and continuous up
to $z =0.$

The space of functions having support in $L_2(U_r(0)\setminus
U_1(0))$ is an invariant subspace for the operator $T_2(r).$ The
kernel of this operator has form
$$
K_{\bf n}(p,q)=\frac{{\bf n}}{\sqrt{3}\pi^2}
\frac{1}{|p|^{\frac{1}{2}} |q|^{\frac{1}{2}}(|p|^2+(p, q)+|q|^2)}.
$$

Let ${\bf T}(r)$ be the integral operator acting on
$L_2(U_r(0)\setminus U_1(0))$ with the kernel $K_2(p,q).$

The following lemma was proven in \cite{Abd-Lak}.

\begin{lem}\label{help result} The equality
$$
\lim\limits_{z\to -0}\frac{d(1, {\bf
T}(z))}{|log|z||}=\frac{\gamma_0}{2\pi}
$$
is satisfied, where $\gamma_0$ is a positive solution of the
equation (\ref{gamma0}).
\end{lem}

Now Theorem \ref{main result} follows from Lemmas \ref{Birm-Schw},
\ref{A0+A1}--\ref{help result}.

\vspace{0.2cm}

{\bf ACKNOWLEDGMENTS.}  The author would like to thank the Abdus
Salam International Centre for Theoretical Physics, Trieste, Italy,
for the kind hospitality and support, and the Commission on
Development and Exchanges of the International Mathematical Union
for the travel grant.

\vspace{0.2cm} {\sc
Samarkand State University,

Department of Physics and Mathematics,

15 University Boulevard, Samarkand, 140104, Uzbekistan}

\vspace{0.1cm} {\it E-mail adress: tulkin$_{-}$rasulov@yahoo.com}


\begin{thebibliography}{9}

\bibitem{Abd-Lak} J. I. Abdullaev and S. N. Lakaev. Asymptotics of the discrete
spectrum of the three-particle Schr\"{o}dinger difference operator
on a lattice. Theor. Math. Phys. \textbf{136} (2003), No. 2,
1096--1109.

\bibitem{AHW} S. Albeverio, R. H\"{o}egh-Krohn and T. T. Wu. A class of exactly
solvable three-body quantum mechanical problems and the universal
low energy behavior. Phys. Lett. A {\bf 83} (1981), No. 3, 105--109.

\bibitem{Alb-Lak-Mum} Sergio Albeverio, Saidakhmat N. Lakaev and Zahriddin I. Muminov.
Schr\"{o}dinger Operators on Lattices. The Efimov Effect and
Discrete Spectrum Asymptotics. Ann. Henri Poincar\'{e}. \textbf{5}
(2004), 743--772.

\bibitem{ALM-3} S. Albeverio, S. N. Lakaev and Z. I. Muminov.
On the Number of Eigenvalues of a Model Operator Associated to a
System of Three-Particles on Lattices. Russian J. Math. Phys. {\bf
14} (2007), No. 4, 377--387.

\bibitem{Amad-Nob} R. D. Amado and J. V. Noble. On Efimov's effect: a new pathology of
three-particle systems. Phys. Lett. B. \textbf{35} (1971), 25--27;
II. Phys. Lett. D. \textbf{5} (1972), No. 3, 1992--2002.

\bibitem{Dell-Fig-Teta} G. F. Dell'Antonio, R. Figari, A. Teta. Hamiltonians for systems
of $N$ particles interacting through point interactions. Ann. Inst.
Henri Poincar\'{e}, Phys. Theor. \textbf{60} (1994), No. 3,
253--290.

\bibitem{Efimov} V. Efimov. Energy levels arising from resonant
two-body forces in a three-body system. Phys. Lett. B \textbf{33}
(1970), No. 8, 563--564.

\bibitem{L-2} S. N. Lakaev. On an infinite number of three-particle bound
states of a system of quantum lattice particles. Theor. Math. Phys.
{\bf 89} (1991), No. 1, 1079--1086.

\bibitem{L-1} S. N. Lakaev. The Efimov Effect of a System of Three
Identical Quantum Lattice Particles. Funct. Anal. Appl. {\bf 27}
(1993), No. 3, 166--175.

\bibitem{L-M} S. N. Lakaev and M. \'{E}. Muminov. Essential and
discrete spectra of the three-particle Schr\"{o}dinger operator on a
lattices. Theor. Math. Phys. \textbf{135} (2003), No. 3, 849--871.

\bibitem{Ovch-Sig} Yu. N. Ovchinnikov and I. M. Sigal. Number of bound states of
three-body systems and Efimov's effect. Ann. Phys. \textbf{123}
(1979), No. 2, 274--295.

\bibitem{Sob} A. V. Sobolev. The Efimov effect. Discrete spectrum asymptotics.
Comm. Math. Phys. \textbf{156} (1993), 101--126.

\bibitem{Tam-1} H. Tamura. The Efimov effect of three-body Schr\"{o}dinger
operators. J. Func. Anal. \textbf{95} (1991), No. 2, 433--459.

\bibitem{Tam-2} H. Tamura.
The Efimov effect of three-body Schr\"{o}dinger operators:
asymptotics for the number of negative eigenvalues. Nagoya Math. J.
{\bf 130} (1993), 55--83.

\bibitem{Wang} X. P. Wang. On the existence of the $N$ - body Efimov effect.
J. Func. Anal. \textbf{209} (2004), 137--161.

\bibitem{Yaf} D. R. Yafaev. On the theory of the discrete spectrum of
the three-particle Schr\"{o}dinger operator. Math. USSR-Sb.
\textbf{23} (1974), 535--559.

\end{thebibliography}
\end{document}